\documentclass[11pt]{amsart}
\usepackage{amssymb}
\newtheorem{theorem}{Theorem}[section]
\newtheorem{lemma}[theorem]{Lemma}
\newtheorem{proposition}[theorem]{Proposition}
\newtheorem{corollary}[theorem]{Corollary}
\theoremstyle{definition}

\theoremstyle{remark}

\numberwithin{equation}{subsection}
\newfont{\ajb}{eufm10 at12pt}
\newfont{\aj}{eufm10 at10pt}
\newfont{\ajk}{eufm10 at8pt}
\newfont{\kh}{msbm10 at10pt}
\newfont{\khk}{msbm10 at 8pt}
\newcommand{\bea}{\begin{eqnarray*}}
\newcommand{\eea}{\end{eqnarray*}}

\begin{document}
\title[A Class of Distal Functions$\cdots$ ]{A Class of distal Functions on Semitopological
Semigroups}
\author{A. Jabbari }
\address{ Department of Mathematics, Ferdowsi University of
Mashhad, P. O. Box 1159, Mashhad 91775,
Iran}\email{shahzadeh@math.um.ac.ir}
\author{H.\ R.\ E.\ Vishki}
\address{Department of Mathematics, Ferdowsi University of Mashhad\\
P. O. Box 1159, Mashhad 91775, Iran; \newline Centre of Excellence
in Analysis on Algebraic Structures (CEAAS), Ferdowsi University of
Mashhad, Iran.} \email{vishki@ferdowsi.um.ac.ir}

\subjclass[2000]{ 54H20, 54C35, 43A60}
\keywords{semitopological semigroup, distal function, strongly almost periodic
function, semigroup compactification, m-admissible subalgebra.}

\date{}

\dedicatory{}

\commby{}

%%%----------------------------------------------------------------------

\begin{abstract}
 The norm closure of the algebra generated
by the set  $\{n\mapsto {\lambda}^{n^k}:$ $\lambda\in{\mathbb {T}}$
and $k\in{\mathbb{N}}\}$ of functions on $({\mathbb {Z}}, +)$ was
studied  in \cite{S} (and was named as the Weyl algebra). In this
paper, by a fruitful result of Namioka, this algebra is generalized
for a general semitopological semigroup and, among other things, it
is shown that the elements of the involved algebra are distal. In
particular, we examine this algebra for $({\mathbb {Z}}, +)$ and
(more generally) for the discrete (additive)  group of any countable
ring. Finally, our results are treated for a bicyclic semigroup.
\end{abstract}

%%%----------------------------------------------------------------------
\maketitle
%%%----------------------------------------------------------------------
\section{Introduction} Distal functions on topological groups were extensively studied by Knapp ~\cite{K}.
The norm closure of the algebra generated by the set $F=\{n\mapsto
{\lambda}^{n^k}$: $\lambda\in{\mathbb {T}}$ and $k\in{\mathbb{N}}\}$
of functions on $({\mathbb {Z}}, +)$ was called the Weyl algebra by
E. Salehi in \cite{S}. Knapp, ~\cite{K}, showed that all of the
elements of $F$ are distal on $({\mathbb {Z}}, +)$. Also Namioka
~\cite[Theorem 3.6]{N} proved the same result by using a very
fruitful result (\cite[Theorem 3.5]{N}) which played an important
role for the construction of this paper. By the above mentioned
results of Knapp and Namioka, all elements of the Weyl algebra are
distal, however it dose not exhaust all distal functions on
$({\mathbb {Z}}, +)$; \cite[Theorem 2.14]{S}. In this paper, we
generalize the notion of Weyl algebra to an arbitrary
semitopological semigroup and also we show that all elements of the
involved algebra are distal. In particular, our method  provides a
convenient way to deduce a result of M. Filali~\cite{Fi} on the
distality of the functions $\chi(q(t))$, where $\chi$ is a character
on the discrete additive group of a (countable) ring $R$ and $q(t)$
is a polynomial with coefficients in $R$.
\section{Preliminaries}
For the background materials and notations we follow Berglund et al. \cite{BJM} as
much as possible. For a semigroup $S$,  the right translation $r_{t}$ and the
left translation $l_{s}$ on $S$ are defined by
$r_{t}(s)=st=l_{s}(t)$, ($s, t\in S$). The semigroup $S$, equipped with a topology, is said to be
right topological if all of the right translations are continuous,
semitopological if all of the left and right translations are
continuous.  If $S$ is a right topological semigroup then the set
$\Lambda(S)=\{s\in S: l_{s}$ is continuous$\}$ is called the
topological centre of $S$.\\
Throughout this paper, unless otherwise stated, $S$ is a
semitopological semigroup. The space of all bounded continuous
complex valued functions on $S$ is denoted by $C(S)$. For
$f\in{{C}}(S)$ and $s\in S$ the right (respectively, left)
translation of $f$ by $s$ is the function $R_{s}f=f\circ r_{s}$
(respectively,
$L_{s}f=f\circ l_{s}$).\\
A left translation invariant unital $C^\ast$-subalgebra $F$ of $C(S)$
(i.e., $L_sf\in F$ for all $s\in S$ and $f\in F$) is called $m$-admissible if the function
$s\mapsto (T_\mu f)(s)=\mu(L_sf)$ belongs to $F$ for all $f\in F$
and $\mu\in S^F$(=the spectrum of F). If $F$ is $m$-admissible then  $S^F$ under the
multiplication $\mu\nu=\mu\circ T_\nu$ ($\mu, \nu\in S^F$),
furnished with the Gelfand topology is a compact Hausdorff right topological semigroup and it  makes $S^F$ a
compactification (called the $F$-compactification) of $S$.\\
Some of the usual $m$-admissible subalgebras of $C(S)$, that are
needed in the sequel, are the left multiplicatively continuous,
strongly almost periodic and distal functions on $S$. These are
denoted by $LMC(S)$, $SAP(S)$ and $D(S)$, respectively. Here and
also for other emerging spaces when there is no risk of confusion,
we have suppressed the letter $S$ from the notation.\\
The interested reader may refer to \cite {BJM} for ample information about these $m$-admissible  subalgebras and the properties of their corresponding compactifications.

\section{Main results}
The idea of defining our new algebras $W_k$ and
$W$, in the form given
below, came from a nice result of Namioka~\cite[Theorem 3.5]{N}.\\
Let $\Sigma=\{T_{\mu}
:LMC(S)\rightarrow LMC(S); \mu\in S^{{LMC}}\}$. Let $F_0$ be the set
of all constant functions of modulus 1. For every $k\in\mathbb{N}$
assume that we have defined $F_i$ for $i=1,2,\ldots,k-1$ and define
$F_k$ by
\[F_{k}=\{f\in
LMC : |f|=1  {\rm \ and \ for \ every \ }  \sigma\in\Sigma,
\sigma(f)=f_{\sigma}f, {\rm \ for \ some \ } f_{\sigma}\in
F_{k-1}\};\] It is clear from definitions that
${{F}}_{k}\subseteq{{F}}_{k+1}$, for all $k\in\mathbb{N}\cup\{0\}$.
Let ${{W}}_{k}$ and $W$ be the norm closure of the algebras
generated by ${{F}}_{k}$ and ${\bigcup}_{k\in\mathbb {N}}^{} F_k$ in
${{C}}(S)$, respectively; then trivially
${{W}}_{k}\subseteq{{W}}_{k+1}\subseteq W$. Hence, $W$ is the
uniform closure of the algebra $\bigcup_{k\in{\mathbb{N}}}W_k$. It
is also readily verified that $W$ is the direct limit of the family
$\{W_i:~i\in{\mathbb{N}}\}$ (ordered by inclusion, and with the
inclusion maps as morphisms). \noindent From now on, we assume that
$k\in\mathbb{N}$ is arbitrary. We leave the following simple
observations without proof.
\begin{proposition}\label{I}

 $(i)$ For every $f\in F_k$ and  $\sigma\in\Sigma$ the function $f_{\sigma}$ with the
property $\sigma(f)=f_{\sigma}f$ is unique.

$(ii)$ For every $f, g\in
{{F}}_{k}$ and $\sigma\in\Sigma$,
$(fg)_{\sigma}=f_{\sigma}g_{\sigma}$. In particular,
${{F}}_{k}$ is a multiplicative subsemigroup of $LMC$.

$(iii)$ $F_k$ is conjugate closed; in other words, if $f\in F_k$ then $\overline{f}\in F_k$.

$(iv)$ $F_k$ contains the constant functions.
\end{proposition}
\begin{lemma}\label{II}  The set $F_k$ is left translation invariant and it is also invariant  under $\Sigma$; in other words, $L_S(F_k)\subseteq F_k$ and $\Sigma(F_k)\subseteq F_k$.
\end{lemma}
\begin{proof} A direct verification reveals that  $F_k$ is left translation invariant. For the invariance under $\Sigma$ let  $f\in {{F}}_{k}$ and  $\sigma\in\Sigma$, the equality ${\sigma}(f)=f_{\sigma}f$ for some $f_\sigma\in F_{k-1}$ implies that $|\sigma(f)|=1$ and so for every $\tau\in\Sigma$, \[\tau(\sigma(f))=\tau(f_\sigma f)=\tau(f_\sigma)\tau(f)=(({f_\sigma})_\tau f_\sigma)(f_\tau f)=({(f_\sigma})_\tau f_\tau)(f_\sigma f)=(({f_\sigma})_\tau f_\tau)\sigma(f).\] Since $({f_\sigma})_\tau f_\tau\in F_{k-1}$ we have $\sigma(f)\in F_k$; in other words $\Sigma(F_k)\subseteq F_k$, as required.
\end{proof}
\begin{lemma}\label{III}  All
elements of ${{F}}_{k}$ remain fixed under the
idempotents of $\Sigma$.
\end{lemma}
\begin{proof} It is easily seen that the result holds for $k=1$. Assume that $k>1$ and that the result holds for $k-1$. Let $f\in
F_{k}$ and let $\varepsilon\in\Sigma$ be an idempotent, then
$\varepsilon(f)=f_{\varepsilon}f$ for some $f_{\varepsilon}\in
F_{k-1}$. Therefore
\[f_{\varepsilon}f=\varepsilon(f)={\varepsilon}^2
(f)={\varepsilon}({\varepsilon}(f))={\varepsilon}(f_{\varepsilon}f)=\varepsilon(f_\varepsilon)\varepsilon(f)=f_\varepsilon(f_\varepsilon
f)={f_{\varepsilon}}^2f;\] hence ${f_{\varepsilon}}=1$ and
${\varepsilon}(f)=f$, as claimed.
\end{proof}
\begin{lemma}\label{IV}
 ${{F}}_{k}\subseteq {{D}}$ .
\end{lemma}
\begin{proof}  Let $f\in
{{F}}_{k}$. To show that $f\in{{D}}$, using \cite[Lemma 4.6.2]{BJM}, it is
enough to show that $\varepsilon\sigma(f)=\sigma(f)$ for each
$\sigma\in\Sigma$ and each idempotent $\varepsilon$ in
$\Sigma$. By Lemma~\ref{II}, $\sigma(f)\in
{{F}}_{k}$, so that  Lemma~\ref{III} implies that
$\varepsilon(\sigma(f))=\sigma(f)$, as required.
\end{proof}
Using parts $(iii)$ and $(iv)$ of Proposition~\ref{I}, $W_k$ and $W$
are unital $C^*-$subalgebras of $C(S)$. The following result shows
that these are indeed $m-$admissible subalgebras of $D$.
\begin{theorem}\label{iV}  For every semitopological semigroup $S$, ${{W}}_{k}$  and ${W}$ are
$m$-admissible subalgebras of $D(S)$.
\end{theorem}
\begin{proof} For the $m-$admissibility of $W_k$ it is enough to show that it is left translation invariant and also invariant under $\Sigma$.
Let $\langle{F_k}\rangle$ be the algebra generated by $F_k$.
Lemma~\ref{II} implies that $L_S(\langle{F_k}\rangle)\subseteq
\langle{F_k}\rangle$ and also $\Sigma(\langle{F_k}\rangle)\subseteq
\langle{F_k}\rangle$. For every $f\in{{W}}_{k}$ there exists a
sequence $\{f_n\}\subseteq \langle{F_k}\rangle$  which converges (in
the norm of $C(S)$) to $f$. Let $\sigma\in\Sigma$ and $s\in S$, then
the inequalities  $\|L_s(f_n)-L_s(f)\|\leq\|f_{n}-f\|$ and
$\|\sigma(f_n)-\sigma(f)\|\leq\|f_{n}-f\|$ imply that
$L_s(f_{n})\rightarrow L_s(f)$ and $\sigma(f_n)\rightarrow
\sigma(f)$, respectively. Since for each $n\in \mathbb {N}$,
$L_s(f_n)$ and $\sigma(f_n)$ lie in $\langle{{{F}}_{k}}\rangle$, we
have $L_s(f)\in W_k$ and also $\sigma(f)\in{{W}}_{k}$. It follows
that ${{W}}_{k}$ is $m$-admissible. A similar argument may apply for
the $m-$admissibility of $W$. The fact that $W_k$ and $W$ are
contained in $D$ follows trivially from Lemma~\ref{IV}.
\end{proof}
The next result gives $S^{{W}}$ in terms of the subdirect product of
the family $\{S^{{{W}}_{k}}: k\in{\mathbb {N}}\}$. For a full
discussion of the subdirect product of compactifications one may
refer to \cite[Section 3.2]{BJM}.

\begin{proposition} The compactification $S^{{W}}$ is the subdirect product of the  family
$\{S^{{{W}}_{k}}: k\in{\mathbb
{N}}\}$; in symbols, $S^{{W}}=\bigvee\{S^{{{W}}_{k}}: k\in{\mathbb {N}}\}$.
\end{proposition}
\begin{proof} The family (of homomorphisms)
$\{\pi_k:S^{{W}}\rightarrow S^{{{W}}_{k}}; k\in{\mathbb {N}}\}$,
where for each $\mu\in S^{{W}}$, $\pi_k(\mu)=\mu |_{{{W}}_{k}}$,
separates the points of $S^{{W}}$, because for given $\mu,\nu\in
S^W$ with $\mu\neq\nu$ (on $W$) one has $\mu\neq\nu$ of
$F=\bigcup_{k\in{\mathbb{N}}}F_k$, hence there exists a natural
number $j$ and an element $f\in F_j$ such that $\mu(f)\neq\nu(f)$.
Therefore $\mu |_{{{W}}_{j}}\neq\nu |_{{{W}}_{j}}$, that is
$\pi_j(\mu)\neq\pi_j(\nu)$. Now the conclusion follows from
\cite[Theorem 3.2.5]{BJM}.
\end{proof}
\begin{proposition}\label{V}
$(i)$ For every  abelian semitopological semigroup $S$,   $SAP(S)\subseteq
W_k(S)$.\\
 \indent $(ii)$ For every abelian semitopological semigroup $S$ with a left
identity, ${{W}}_1(S)={SAP(S)}$.
\end{proposition}
\begin{proof} $(i)$ Since $S$ is abelian ${{SAP}}(S)$ is the
closed linear span of the set of all continuous characters on $S$; see \cite[Corollary 4.3.8]{BJM}.
Let $f$ be any continuous character on $S$ and let $\sigma\in\Sigma$. Then there exists a net $\{s_\alpha\}$ in $S$ such that $\sigma(f)=\lim_\alpha R_{s_\alpha}f$. By passing to a subnet, if necessary, we may assume that $f(s_\alpha)$ converges to some $\lambda_\sigma\in\mathbb{T}$. Therefore for every $s\in S$, $\sigma(f)(s)=\lim_\alpha R_{s_\alpha}f(s)=\lim_\alpha f(ss_\alpha)=f(s)\lim_\alpha f(s_\alpha)=f(s)\lambda_\sigma$. Hence $\sigma(f)=\lambda_\sigma f$ or equivalently $f\in F_1$.

$(ii)$ By part $(i)$ it is enough to show that ${{W}}_1\subset{{SAP}}$. Indeed we are going to show that $F_1\subset{{SAP}}$; for this end, let $f\in{{F}}_1$ and let $s\in S$, then
$R_sf=\lambda_s f$ for some $\lambda_s$ in ${\mathbb{T}}$. Let $e$
be a left identity of $S$, then for each $s$ in $S$,
$f(s)=R_sf(e)=\lambda_s f(e)$. Let $h=\frac{1}{f(e)}f$, then $h$ is
a continuous character on $S$. But $f=f(e)h$, now using the fact
that ${{SAP}}$ is the closed linear span of continuous characters
of $S$ we have $f\in{{SAP}}$, as required.
\end{proof}
As a consequence of the latter result we have

\begin{corollary} For any compact abelian
topological group $G$, $W_k(G)=W(G)={{C}}(G)$.
\end{corollary}
\begin{proof} Since for every compact topological group
$G$, ${{SAP}}(G)={{C}}(G)$, \cite[Theorem 4.3.5]{BJM}, the result follows from
the last proposition.
\end{proof}

\section{Examples}

{\bf  Example  (a).} Here we restrict our discussion to the
discrete group $(\mathbb{Z}, +)$ and examine $W$ and $W_k$ for this
particular case, which were studied extensively by  Salehi in
~\cite{S}. Note that although we would deal with countable discrete
rings in part (b), but the proofs on ${\mathbb{Z}}$ are more
interesting and characterizations of $F_k({\mathbb{Z}})$ are more
explicit. We commence with the next key lemma which characterizes
$F_k$ in $l^\infty(\mathbb{Z})$.
\begin{lemma} The set ${{F}}_{k}({\mathbb {Z}}, +)$ is the
(multiplicative) sub-semigroup of $l^{\infty}({\mathbb {Z}})$
generated by the set $\{n\mapsto {\lambda}^{n^i}:$
$\lambda\in{\mathbb {T}}$,
$i=0, 1,..., k\}$.
\end{lemma}
\begin{proof} For each $k\in{\mathbb{N}}$, let
$A_k$ denote the multiplicative sub-semigroup of
$l^{\infty}({\mathbb {Z}})$ generated by the set $\{n\mapsto
{\lambda}^{n^i}$, $\lambda\in{\mathbb {T}}$ and $i=0, 1,..., k\}$.
For $k=1$ a direct verification reveals that $A_1\subseteq F_1$; for
the reverse inclusion let $f\in F_1$. Then for some
$\lambda\in{\mathbb {T}}$, $R_{1}f={\lambda}f$, hence
$f(1)=R_{1}f(0)={\lambda}f(0)={\lambda}\lambda_1$, in which
$\lambda_1=f(0)$. Also
$f(2)=R_{1}f(1)={\lambda}f(1)={\lambda}^{2}\lambda_1$, by induction
it is easily proved that for each $n\in{\mathbb {N}}$,
$f(n)={\lambda}^{n}\lambda_1$. Let $R_{-1}f={\beta}f$, where
$\beta\in{\mathbb{T}}$, then $f(-1)=R_{-1}f(0)={\beta}f(0)$. But
$f(1)=R_{-1}f(2)={\beta}f(2)={\beta}{\lambda}^{2}\lambda_1$,
therefore ${\lambda}\lambda_1={\beta}{\lambda}^{2}\lambda_1$, hence
$\beta={\lambda}^{-1}$ and for all $n\in{\mathbb {Z}}$,
$f(n)={\lambda}^{n}\lambda_1$. Thus
$f\in A_1$, and so $F_1=A_1$.\\
Let $k\geq 2$ and assume that $A_{k-1}=F_{k-1}$. Let $n\in{\mathbb
{Z}}$ and $\lambda\in{\mathbb {T}}$ and let $f\in A_k$ and assume
(without loss of generality) that $f(n)={\lambda}^{n^{k}}$, then for
given $\sigma={\lim}_{\alpha}m_{\alpha}\in\Sigma$ we have
$\sigma(f)(n)={\lim}_{\alpha}R_{m_{\alpha}}f(n)={\lim}_{\alpha}{\lambda}^{(n+m_{\alpha})^{k}}=
f(n)f_{\sigma}(n)$, in which
$f_{\sigma}(n)={\mu}_{1}^{n^{k-1}}{\mu}_{2}^{n^{k-2}}...{\mu}_{k-1}^{n}{\mu}_{k}$,
where (by going through a sub-net of $\{m_\alpha\}$, if necessary)
${\mu}_{i}={\lim}_{\alpha}{\lambda}^{({}_{i}^{k}){m_{\alpha}}^{i}}$,
for $i=1, 2,..., k$. But $f_{\sigma}\in A_{k-1}=F_{k-1}$, so $f\in
F_k$. Hence $A_k\subseteq F_k$. Now  let $f\in F_k$, we have to show
that $f\in A_k$. We have $R_{1}f=f_{1}f$, for some $f_{1}\in
F_{k-1}=A_{k-1}$. Since $f_1\in A_{k-1}$ we may assume that
$f_{1}(n)={\lambda}_{1}^{n^{k-1}}{\lambda}_{2}^{n^{k-2}}...{\lambda}_{k-1}^{n}{\lambda}_{k}$,
where ${\lambda}_{1}, {\lambda}_{2},..., {\lambda}_{k}\in{\mathbb
{T}}$. Then $f(1)=R_{1}f(0)=f_{1}(0)f(0)$ and
$f(2)=R_{1}f(1)=f_{1}(1)f(1)=f_{1}(1)f_{1}(0)f(0)$, and by
induction, \bea
f(n)&=&({\lambda}_{1}^{(n-1)^{k-1}}{\lambda}_{2}^{(n-1)^{k-2}}...\lambda_{k-1}^{n-1}{\lambda}_{k})
({\lambda}_{1}^{(n-2)^{k-1}}{\lambda}_{2}^{(n-2)^{k-2}}...{\lambda}_{k-1}^{n-2}{\lambda}_{k})
...({\lambda}_{1}\lambda_2...{\lambda}_{k})
({\lambda}_{k})f(0)\\
&=&{\lambda}_{1}^{{\sum}_{j=1}^{n-1}j^{k-1}}{\lambda}_{2}^{{\sum}_{j=1}^{n-1}j^{k-2}}
...{\lambda}_{k-1}^{{\sum}_{j=1}^{n-1}j}{\lambda}_{k}^{n}f(0). \eea
So for each $i=0,1,2,\ldots,k-1$ the power of $\lambda_{k-i}$ is a
polynomial in $n$ of degree $i+1$. Hence the power of $\lambda_1$
(which has the maximum degree) is a polynomial of degree $k$. It
follows that $f\in A_k$ and the proof is complete by induction.
\end{proof}
As an immediate consequence of the latter lemma we have the next
theorem which characterizes our algebras for the additive group
$\mathbb Z$.
\begin{theorem}
$(i)$  $W_{k}({\mathbb {Z}}, +)$ coincides with the norm closure of
the algebra generated by the set of functions $\{n\mapsto
{\lambda}^{n^i}:$ $\lambda\in{\mathbb {T}}$, $i=0, 1,..., k\}$.

$(ii)$  $W({\mathbb {Z}}, +)$ coincides with the norm closure of the
algebra generated by the set of functions
$\{n\mapsto{\lambda}^{n^k}: \lambda\in{\mathbb {T}}$, and
$k\in{\mathbb {N}}\}$, that is, $W({\mathbb {Z}}, +)$ coincides with
the Weyl algebra.
\end{theorem}

{\bf Example (b).} Let $R$ be a countable discrete ring.
Let $\chi$ be an arbitrary character on the additive group
$(R_d,+)$, where $R_d$ denotes $R$ with the discrete topology.
Without loss of generality, assume that $R$ is abelian. We are going
to show that for each $s$ in $R$ the function $f(t)=\chi(st^3)$
belongs to $F_3(R_d,+)$. To this end, let $\sigma\in\Sigma$. Thus
there exists a sequence $s_n$ in $R$ such that for each $h\in
l^\infty(R_d)$, $\sigma(h)(t)=\lim_n
R_{s_n}h(t)=\lim_n h(t+s_n)$.\\
Thus $\sigma f(t)=\lim_n
f(t+s_n)=\lim_n\chi(s{(t+s_n)}^3)=f(t)f_\sigma(t)$, in which
$f_\sigma(t)=\lim_n\chi(ss_n^3+3ss_n t^2+3ss_n^2 t)$. Since $R$ is
countable, by the diagonal process, there exists a subsequence, say
$s_n$, of the sequence $s_n$ such that, for all $t$ in $R$, all of
the limits $\lim_n\chi(ss_n^3)$, $\lim_n\chi(ss_n t^2)$ and
$\lim_n\chi(ss_n^2 t)$ exist. (In fact, one may first choose a
subsequence of $s_n$, if necessary, such that $\lim_n\chi(ss_n^3)$
exist. Let $R=\{x_1,x_2,x_3,\ldots\}$. Choose a subsequence of
$s_n$, say $s_{1,n}$ such that both limits $\lim_n\chi(ss_{1,n}
x_1^2)$ and $\lim_n\chi(ss_{1,n}^2 x_1)$ exist. This time, choose a
subsequence of $s_{1,n}$, say $s_{2,n}$, such that both limits
$\lim_n\chi(ss_{2,n} x_2^2)$ and $\lim_n\chi(ss_{2,n}^2 x_2)$ exist.
Continue this process and choose the resulting sequence $s_{n,n}$ on
the diagonal, which is eventually our desired subsequence, say
$s_n$). Hence for each $t\in R$,
\[f_\sigma(t)=\lim_n\chi(ss_n^3)\lim_n\chi^3(ss_n^2
t)\lim_n\chi^3(ss_n t^2).\] By definition, to prove that $f\in F_3$
it is enough to show that $f_\sigma\in F_2$. To see this, let
$\tau\in\Sigma$ be arbitrary. Then there exists a sequence $t_m$ in
$R$ such that for each $h\in l^\infty(R_d)$, $\tau(h)(t)=\lim_m
R_{t_m}h(t)=\lim_m h(t+t_m)$.\\
Thus $\tau f_\sigma(t)=\lim_m
f_\sigma(t+t_m)=f_\sigma(t)f_{\sigma\tau}(t)$, where
\[f_{\sigma\tau}(t)={(f_\sigma)_\tau}(t)=\lim_m\lim_n\chi^3(ss_n^2
t_m)\lim_m\lim_n\chi^3(ss_n t_m^2+2ss_n t_m t).\] $R$ is countable,
hence by going through a subsequence of $t_m$ (by using the diagonal
process) one may assume that for all $t$ in $R$ the limits
$\lim_m\lim_n\chi(ss_n t_m^2)$ and $\lim_m\lim_n\chi(ss_n t_m t)$
exist. Therefore
\[f_{\sigma\tau}(t)=\lim_m\lim_n\chi^3(ss_n^2
t_m)\lim_m\lim_n\chi^3(ss_n t_m^2)\lim_m\lim_n\chi^6(ss_n t_m t).\]
Again by definition, to prove that $f_\sigma\in F_2$ it is enough to
show that $f_{\sigma\tau}\in F_1$. Let $\xi=\lim_l u_l\in\Sigma$,
then it follows from the above equation that
\[\xi(f_{\sigma\tau})(t)=f_{\sigma\tau}(t)\lim_l\lim_m\lim_n\chi^6(ss_n
t_m u_l).\] That is, $\xi(f_{\sigma\tau})=\lambda f_{\sigma\tau}$
where $\lambda=\lim_l\lim_m\lim_n\chi^6(ss_n t_m u_l)\in
F_0={\mathbb{T}}$. Hence $f_{\sigma\tau}\in F_1$ and so $f_\sigma\in
F_2$ and this implies that $f\in F_3$. Our claim is now established.
By using the above method, one may prove that for each
$k\in{\mathbb{N}}$ and $s\in R$ the function
$t\rightarrow\chi(st^k)$ is an element of $F_k$.\\
Briefly speaking, the above argument and Lemma ~\ref{IV} imply that:
\begin{corollary}
If $R$ is a countable discrete ring, then for each character $\chi$
on the discrete additive group of $R$ the function $\chi(q(t))$, in
which $q(t)$ is a polynomial with coefficients in $R$, belongs to
$W(R_d,+)$ and is also a distal function. \end{corollary}

It should be remarked that the distality of the functions $\chi(q(t))$ was
first proved by Filali~\cite{Fi}
without the countability condition on $R$.\\

{\bf Example (c).} $(i)$ If $S$ contains a right zero element, i.e.
there exists $t\in S$ such that $st=t$ for all $s\in S$, then for
$f\in F_1$ there exists $\lambda_t\in{\mathbb{T}}$ such that
$R_tf=\lambda_tf$, hence for all $s\in S$,
$f(t)=f(st)=R_tf(s)=\lambda_tf(s)$, that is
$f=f(t)/\lambda_t\in{\mathbb{T}}$. Therefore $F_1={\mathbb{T}}$ and
so $F_k={\mathbb{T}}$ for all $k\in{\mathbb{N}}$. It follows that
for such a semigroup $S$, $W_k(S)=W(S)=$ the set of constant
functions.

$(ii)$ If $S$ is a left zero semigroup (i.e. $st=s$ for all $s,t\in
S$), then for each function $f\in LMC(S)$ we have $\sigma(f)=f$ for
all $\sigma$ in $\Sigma$, and so if $|f|=1$, then $f\in F_1$. Hence
for all $k\in {\mathbb{N}}$, $W=W_k=W_1=LMC$.\\

Now we examine some of the newly defined algebras on a non-trivial
non-group
semigroup.\\
{\bf Example (d).} Let $S$ be the bicyclic semigroup of
\cite[Example 2.10]{BJM}, i.e. $S$ is a semigroup generated by
elements $1$, $p$ and $q$, where $1$ is the identity and $p$ and $q$
satisfy $pq=1\neq qp$. The relation $pq=1$ implies that any member
of $S$ may be written uniquely in the form $q^m p^n$, where
$m,n\in{\mathbb{Z}}^+$ and $p^0=q^0=1$. \\
We are going to show that
\[F_1(S)=\{q^mp^n\mapsto\mu^{r}\nu^{1-r}:~~r=m-n~and
~\mu,~\nu\in{\mathbb{T}}\}\] To see this, let $f\in F_1$, then for
each $s\in S$, $R_sf=\lambda_sf$ for some
$\lambda_s\in{\mathbb{T}}$. Hence $f(q)=\lambda_qf(1)$,
$f(p)=\lambda_pf(1)$ and $f(1)=f(pq)=R_qf(p)=\lambda_qf(p)$,
therefore $\lambda_p\lambda_q=1$. It is also readily seen that
$\lambda_{q^mp^n}=\lambda_q^{m-n}$. One may use induction to simply
prove that for each $n\in{\mathbb{Z}}^+$,
$f(qp^n)=f(q)(\frac{f(p)}{f(1)})^n$, and then use the latter to show
(again  by induction on $m$) that
$f(q^mp^n)=f(p)^{n-m}f(1)^{1-(n-m)}$. But $f(p)f(q)=f(1)^2$, so
$f(q^mp^n)=f(q)^{m-n}f(1)^{1-(m-n)}$. Hence it is enough to take
$\mu=f(q)$ and $\nu=f(1)$. The converse inclusion is easily
verified.\\

To prove the next theorem, the following lemma is needed.
\begin{lemma}\label{l:5}
Let $S$ be the bicyclic semigroup described above. If $f\in F_1(S)$,
then $f(p)f(q)=f(1)^2$.
\end{lemma}
\begin{proof}
Let $f\in F_1(S)$, then there exist constants $\lambda_p$ and
$\lambda_q\in{\mathbb{T}}$ such that $R_pf=\lambda_pf$ and
$R_qf=\lambda_qf$. Hence $f(p)f(q)=\lambda_p\lambda_qf(1)^2$. Thus
it is enough to show that $\lambda_p\lambda_q=1$. But,
$f(1)=f(pq)=R_qf(p)=\lambda_qf(p)=\lambda_q\lambda_pf(1)$, that is
$\lambda_q\lambda_p=1$, and the result follows.
\end{proof}
\begin{theorem}\label{t:5} Let $S$ be the bicyclic semigroup generated by $1, p$ and $q$, where $1$ is the identity and $pq=1\neq qp$.
Then $W_1(S)$ and $W_2(S)$ are the norm closure of the algebras
generated by the sets \[\{q^mp^n\mapsto\mu^{r}\nu^{1-r}:~~r=m-n~and
~\mu,~\nu\in{\mathbb{T}}\}\] and
\[\{q^mp^n\mapsto\lambda^{\frac{r^2-r}{2}}\mu^{\frac{r^2+r}{2}}\nu^{1-r^2},~~r=m-n~
and~\lambda,~\mu,~\nu\in{\mathbb{T}}\}\] respectively.
\end{theorem}
\begin{proof} By what we already discussed, the result is clear for
$W_1(S)$. To complete the proof, it is enough to show that
\[F_2(S)=\{q^mp^n\mapsto\lambda^{\frac{r^2-r}{2}}\mu^{\frac{r^2+r}{2}}\nu^{1-r^2},~~r=m-n~
and~\lambda,~\mu,~\nu\in{\mathbb{T}}\}\] Let $A$ denote the right
hand side of the above equation and let $f\in A$. Then there exist
$\lambda,\mu,\nu\in{\mathbb{T}}$ such that for all
$m,n\in{\mathbb{Z}}^+\cup\{0\}$,
$f(q^mp^n)=\lambda^{\frac{r^2-r}{2}}\mu^{\frac{r^2+r}{2}}\nu^{1-r^2}$
with $r=m-n$. By choosing suitable $m,n\in{\mathbb{Z}}^+\cup\{0\}$
we derive that $\lambda=f(p)$, $\mu=f(q)$ and $\nu=f(1)$. To prove
$f\in F_2(S)$ is to prove that there exist elements $f_p$ and $f_q$
in $F_1(S)$ such that $R_pf=f_pf$ and $R_qf=f_qf$. In fact it is
enough to take
\[f_p(q^mp^n)=[f(q)^{-1}f(1)]^r[f(p)f(1)^{-1}]^{1-r}\] and
\[f_q(q^mp^n)=[f(p)f(q)^2f(1)^{-3}]^r[f(q)f(1)^{-1}]^{1-r}\]
where $r=m-n$. Our discussion preceding the theorem reveals that
both $f_p$ and $f_q$ are elements of $F_1(S)$, therefore $f\in
F_2(S)$.\\
Conversely, let $f\in F_2(S)$. To show $f\in A$ is to show that for
all $m,n\in{\mathbb{Z}}^+\cup\{0\}$,
\[f(q^mp^n)=f(p)^{\frac{r^2-r}{2}}f(q)^{\frac{r^2+r}{2}}f(1)^{1-r^2},~where~r=m-n.~ \ \ \ \ \ (*)\]
 Since $f\in F_2$, there exists
$f_q\in F_1(S)$ such that $R_qf=f_qf$. Therefore, (from the above
Lemma) \[f_q(p)f_q(q)=f_q(1)^2.~ \ \ \ \ \ (I)\] On the other hand,
$f_q(1)=f(q)f(1)^{-1}$, and also $f(1)=f(pq)=R_qf(p)=f_q(p)f(p)$,
thus $f_q(p)=f(p)^{-1}f(1)$. Hence it follows from $(I)$ that
\[f_q(q)=f(p)f(q)^2f(1)^{-3}.~ \ \ \ \ \ (II)\] Fix $n\in{\mathbb{N}}$, then
by using induction on $m$ and the fact that
$f(q^{m+1}p^n)=R_qf(q^mp^n)=f_q(q^mp^n)f(q^mp^n)=f_q(q)^{m-n}f_q(1)^{1-(m-n)}f(q^mp^n)$,
we deduce from $(II)$ and $(*)$ that
\[f(q^{m+1}p^n)=f(p)^{\frac{s^2-s}{2}}f(q)^{\frac{s^2+s}{2}}f(1)^{1-s^2}, ~where~ s=(m+1)-n.\]
The theorem is now established by induction.
\end{proof}
\begin{corollary}
Let $S$ be the bicyclic semigroup generated by $1, p$ and $q$, where
$1$ is the identity and $pq=1\neq qp$. If either $p$ or $q$ is an
idempotent, then $W(S)=W_i(S)={\mathbb{C}}$, for all $i$.
\end{corollary}
\begin{proof}
It is enough to show that $F_1(S)={\mathbb{T}}$. To this end, let
$f\in F_1(S)$, and assume that $p^2=p$. Then there exist
$\lambda_p\in{\mathbb{T}}$ such that $R_pf=\lambda_pf$, hence
$f(p)=f(p^2)=R_pf(p)=\lambda_pf(p)$ and so $\lambda_p=1$. That is,
$R_pf=f$ and $f(p)=f(1)$. It follows from Lemma~\ref{l:5} that
$f(q)=f(1)$. Now the first part of the above theorem implies that
for all $m,n\in{\mathbb{Z}}^+\cup\{0\}$,
$f(q^mp^n)=f(q)^{m-n}f(1)^{1-(m-n)}=f(1)$. Hence
$f=f(1)\in{\mathbb{T}}$. The proof for the case where $q$ is an
idempotent is similar.
\end{proof}
{\bf Remarks.} $(i)$ By the results ~\ref{iV} and ~\ref{V}, for
every abelian semitopological semigroup $S$, $W_k$  and also $W$ lie
between $SAP$ and $D$. It would be desirable to study the structure
of the (right topological abelian group) compactifications $S^{W_k}$
and $S^W$. In particular, it would be more desirable if one could
investigate the size of the topological centres of $S^{W_k}$ and
$S^W$. The latter problem is of particular interest among some
authors, (the interested reader is referred to~\cite{J} and
~\cite{DL}). For an elegant characterization of the topological
centre of the largest compactification of a locally compact group
one may refer to ~\cite{LP} and also~\cite{NU}.

$(ii)$ In \cite[Theorem 2.13]{S} and~\cite[Corollary 3.3.3]{J}, (by
using different methods) it is proved that all elements of the Weyl
algebra $W(\mathbb{Z}, +)$ are uniquely ergodic. One may seek the
same result for $W(S)$, where $S$ is an arbitrary semitopological
semigroup.

$(iii)$ It would be quite interesting to find a general formula for
$F_k(S)$ in Theorem~\ref{t:5}.

\section*{acknowledgment}
The first  author would like to thank Professor Anthony Lau,  at the
University of Alberta, for his encouragement and support through his
NSERC grant A7679.  And the authors would like to thank the kind
referee for the helpful suggestions.

\end{document}